\documentclass[12pt]{amsart}
\usepackage{bm}
 \usepackage{graphicx}

 \newtheorem{thm}{Theorem}[section]

 \theoremstyle{definition}
 
 \theoremstyle{remark}
 
 \numberwithin{equation}{section}

 \newcommand{\tr}{\textbf{tr}}
 \newcommand{\dis}{\textbf{d}}
\newcommand{\ntorus}{\mathbb T^n}

\begin{document}

\title[The Harnack inequality and the JKO scheme]{The Harnack inequality and the Jordan-Kinderlehrer-Otto scheme}

\author{Paul W.Y. Lee}
\email{wylee@math.cuhk.edu.hk}
\address{Room 216, Lady Shaw Building, The Chinese University of Hong Kong, Shatin, Hong Kong}

\date{\today}

\begin{abstract}
We establish a version of the Harnack inequality for the Jordan-Kinderlehrer-Otto scheme of the heat equation on the flat torus. 
\end{abstract}

\maketitle

\section{Introduction}

The Harnack inequality is one the most important inequalities for elliptic and parabolic equations. In the linear elliptic case, the inequality was first proved in \cite{Se}. It was extended to more general elliptic and parabolic equations in  \cite{Mo1, Mo2} (see also \cite{Na, FaSt}). In this paper, we prove a version of the Harnack inequality for the Jordan-Kinderlehrer-Otto (JKO) scheme of the heat equation on the flat torus. 

The JKO scheme is a time discrete scheme introduced in \cite{JoKiOt} using the theory of optimal transportation. It produces time discrete approximations to various gradient flows on the Wasserstein space, the space of probability measures. For the infinite dimensional geometry behind this gradient flow, see \cite{Ot, AmGiSa}. Since its discovery, the JKO scheme has been an intensive area of research. Here is a list which is far from exhaustive \cite{AmGiSa, AmSe, BlCa, CaGa, DeMeSaVe, FiGaYo} (see also \cite{Vi1} for an overview). 

In what follows, we will discuss the detail construction of this scheme in the case of the heat equation in order to give a precise statement of the main result. 

Let $\mu_0$ and $\mu_1$ be two Borel probability measures on the flat $\ntorus$ and let $d$ be the distance function on $\ntorus$. The theory of optimal transportation begins with the following problem which searches for a map $\varphi$ transporting the first mass $\mu_0$ to the second one $\mu_1$ with the minimal amount of total cost: 
\[
\inf_{\varphi_*\mu_0=\mu_1}\int_Md^2(x,\varphi(x))d\mu_0(x)
\]
where the infimum is taken over all Borel maps $\varphi:\ntorus\to \ntorus$ pushing $\mu_0$ forward to $\mu_1$.

The above problem defines a distance function $\dis$, called $L^2$ Wasserstein distance, on the Wasserstein space. More precisely,  
\begin{equation}\label{OT}
\dis^2(\mu_0,\mu_1)=\inf_{\varphi_*\mu_0=\mu_1}\int_Md^2(x,\varphi(x))d\mu_0(x). 
\end{equation}

Assuming the first measure is absolutely continuous with respect to the Lebesgue measure, the existence and uniqueness of minimizer to the above problem (\ref{OT}) was established in \cite{Br,Co, Mc}. This minimizer is called the optimal map that pushes $\mu_0$ forward to $\mu_1$.

For each fixed number $K>0$, each positive integer $N$, and each probability measure $\rho_0\,dx^n$, we consider the following family of minimization problems parametrized by $k=1,...,N$:
\begin{equation}\label{JKOmin}
\inf\left[\frac{1}{2}\dis^2(\rho_{k-1}^N\, dx^n,\rho\, dx^n)+\frac{K}{N}\int_{\ntorus}\rho\log \rho\, dx^n\right],
\end{equation}
where $\rho_0^N=\rho_0$ and the infimum is taken over the set of $L^1$ functions $\rho:\ntorus\to [0,\infty)$ satisfying $\int\rho d\mu=1$.

Assuming that the function $\rho_0$ is in $L^1$, the existence and uniqueness of minimizer to the problem (\ref{JKOmin}) was established in \cite{JoKiOt}. Moreover, if $K$ and $\rho_0$ are fixed, the initial condition $\rho_0$ is in $C^{2,\alpha}$, and $N$ is large enough, then it was shown in \cite{Le} that each $\rho_k^N$ is in $C^{4,\alpha}$.  
The map $\varphi_k^N$ defined by $\varphi_k^N(x)=x+\frac{K}{N}\nabla \log\rho_k^N(x)$ is an optimal map pushing $\rho_{k-1}^N\, dx^n$ forward to $\rho_k^N\,dx^n$. It also satisfies the following Monge-Ampere type equation 
\begin{equation}\label{MA}
\begin{split}
\rho_k^N&=\rho_{k-1}^N(\varphi_k^N)\det\left(d\varphi_k^N\right). 
\end{split}
\end{equation}

The above minimization problem (\ref{JKOmin}) gives a sequence of functions $\Gamma:=\{\rho_k^N|k=0,1,...\}$. Let $u_t^N:\left[0,K\right]\times \ntorus\to [0,\infty)$ be the function defined by
\[
u_t^{N}=\rho_{k}^N
\]
if $t$ is in $\left[\frac{kK}{N},\frac{(k+1)K}{N}\right)$ and $k=0,...,N-1$.

It was shown in \cite{JoKiOt} that $u_t^N$ converges in $L^1$ to the solution of the heat equation with the initial condition $u_0=\rho_0$ as $N$ goes to infinity. This discretization scheme is the so-called Jordan-Kinderlehrer-Otto (JKO) scheme. It was also shown in \cite{Le} that there is a uniform $C^1$ bound for the family $\Gamma$  and so the convergence is improved to $C^{0,\alpha}$ in space (the $C^0$ bound was also observed in \cite{Sa}). In this paper, we establish a uniform lower bound for the second derivatives of the elements in $\Gamma$. More precisely, 

\begin{thm}\label{diffHarn}
Assume that $\rho_0$ is in $C^{2,\alpha}$. There is a constant $1\geq C\geq \frac{1}{2}$ and, for each fixed $K>0$ and each $\rho_0$, there is an integer $N_0>0$ such that 
\begin{equation}\label{diffIneq}
\nabla^2\log u_t^N\geq-\frac{C}{t}I
\end{equation}
for all $N\geq N_0$. 
\end{thm}

By letting $N$ goes to $\infty$ in (\ref{diffIneq}), we recover the matrix differential Harnack inequality proved in \cite{Ha} for the heat equation. This matrix differential Harnack inequality is the matrix analogue of a scalar version proved in \cite{LiYa}. They are called differential Harnack inequalities because one can recover the Harnack inequality by integrating the differential ones along geodesics. 

Theorem \ref{diffHarn} shows that the matrix differential Harnack inequality holds at the level of the JKO scheme. The following is the Harnack inequality for the JKO scheme mentioned at the beginning of the introduction. 

\begin{thm}\label{Harn}
For each fixed $K>0$ and each $\rho_0$, there is an integer $N_0>0$ such that 
\[
\begin{split}
&u_{t_1}^N(x)\leq \left(\frac{t_2+\frac{K}{N}}{t_1}\right)^n\exp\left(\frac{d^2(x, y)}{2\left(t_2-t_1-\frac{K}{N}\right)}\right)u_{t_2}^N(y)
\end{split}
\]
for all $N\geq N_0$. 
\end{thm}

By letting $N$ goes to $\infty$, we recover a version of the Harnack inequality for the heat equation. 

The rest of the paper is devoted to the proof of Theorem \ref{diffHarn} and \ref{Harn}. 

\smallskip

\section{Proof of Theorem \ref{diffHarn}}

Let $\frac{K}{N}f_k^N$ be the $c$-transform of $-\frac{K}{N}\log\rho_k^N$ defined by 
\[
\frac{K}{N}f_k^N(x)=\inf_{y\in\ntorus}\left(\frac{1}{2}d^2(x,y)+\frac{K}{N}\log\rho_k^N(y)\right). 
\]

By definition, the function $f_k^N$ is locally semi-concave. In particular, it is differentiable Lebesgue almost everywhere and the map $\phi_k^N(x):=x-\frac{K}{N}\nabla f_k^N(x)$ is defined Lebesgue almost everywhere. 
By \cite[Lemma 4.1]{Le}, $\rho_k^N$ is bounded uniformly away from $0$ and $\infty$. By \cite[Theorem 9]{Mc}, 
\begin{equation}\label{equal}
\frac{K}{N}f_k^N(x)=\frac{1}{2}d^2(x,\phi_k^N(x))+\frac{K}{N}\log\rho_k^N(\phi_k^N(x))
\end{equation}
for Lebesgue almost all $x$ and 
\begin{equation}\label{ineq}
\frac{K}{N}f_k^N(x)\leq \frac{1}{2}d^2(x, y)+\frac{K}{N}\log\rho_k^N(y)
\end{equation}
for all $x$ and $y$ in $\ntorus$. 

Therefore, for Lebesgue almost all $x$, the function $y\mapsto \frac{1}{2}d^2(x,y)+\frac{K}{N}\log\rho_k^N(y)$ achieves its minimum at $\phi_k^N(x)$. Therefore, 
\begin{equation}\label{sec}
d\varphi_k^N=I+\frac{K}{N}\nabla^2\log\rho_k^N\geq 0
\end{equation}
at $\phi_k^N(x)$. 

By \cite[Theorem 11]{Mc}, 
\begin{equation}\label{inverse}
\phi_k^N(\varphi_k^N(x))=x
\end{equation} 
Lebesgue almost everywhere. Therefore, (\ref{sec}) holds Lebesgue almost everywhere. Since $\varphi_k^N$ is continuous, (\ref{sec}) holds everywhere. 

By combining (\ref{equal}), (\ref{ineq}), and (\ref{inverse}), it follows that 
\begin{equation}\label{pos}
d\varphi_k^N=I+\frac{K}{N}\nabla^2\log\rho_k^N\geq 0. 
\end{equation}
Therefore, by (\ref{MA}) and the uniform bound of $\rho_k^N$, $d\varphi_k^N\geq cI>0$, where $c$ is independent of $k$ and $N$. It also follows that $x\mapsto x+\frac{tK}{N}\nabla\log\rho_k^N(x)=:\psi_k^N(t,x)$ is a diffeomorphism for each $t$ in $[0,1]$. Hence, by the method of characteristics (see \cite{Ev}), the Hamilton-Jacobi equation 
\[
\dot g_k^N+\frac{1}{2}|\nabla g_k^N|^2=0
\]
with initial condition $g_k^N(0,x)=\frac{K}{N}\log\rho_k^N(x)$ has a smooth solution. Moreover, $\psi_k^N$ is the flow of the vector field $\nabla g_k^N$. It follows that 
\begin{equation}\label{HJ}
\begin{split}
\frac{d}{dt}g_k^N(\psi_k^N)&=-\frac{1}{2}|\nabla g_k^N|^2_{\psi_k^N}+\left<\nabla g_k^N,\dot\psi_k^N\right>\\
&=\frac{1}{2}|\nabla g_k^N|^2_{\psi_k^N}. 
\end{split}
\end{equation}

Since $t\mapsto \psi_k^N(t,x)$ is length minimising between its endpoints $x$ and $\varphi_k^N(x)$. It follows from integrating (\ref{HJ}) that 
\begin{equation}\label{g=f}
g_k^N(1,\varphi_k^N(x))=\frac{K}{N}\log\rho_k^N(x)+\frac{1}{2}d^2(x,\varphi_k^N(x))=\frac{K}{N}f_k^N(\varphi_k^N(x)). 
\end{equation}
The last equality follows from (\ref{equal}) and (\ref{inverse}). 

On the other hand, since $\dot\psi_k^N=\nabla g_k^N(\psi_k^N)$, it follows that 
\[
\frac{d}{dt}d\psi_k^N=\nabla^2 g_k^N(\psi_k^N)\, d\psi_k^N. 
\]
Therefore, by (\ref{g=f}), 
\begin{equation}\label{fest}
\begin{split}
\frac{K}{N}\nabla^2f_k^N(\varphi_k^N)&=\nabla^2 g_k^N(\psi_k^N)\Big|_{t=1} \\
&=\frac{K}{N}\nabla^2\log\rho_k^N\left(I+\frac{K}{N}\nabla^2\log\rho_k^N\right)^{-1}. 
\end{split}
\end{equation}

By \cite[Theorem 11]{Mc}, $\phi_k^N$ is the optimal map pushing forward $\rho_{k-1}^N\,dx^n$ to $\rho_{k}^N\,dx^n$ and so it satisfies 
\[
\rho_{k}^N(\phi_k^N)\det(d\phi_k^N)=\rho_{k-1}^N. 
\]

By combining this with (\ref{equal}), we obtain 
\[
f_k^N-\log\rho_{k-1}^N+\log\det\left(I-\frac{K}{N}\nabla^2f_k^N\right) =\frac{N}{2K}d^2(x,\phi_k^N(x))=\frac{K}{2N}|\nabla f_k^N|^2. 
\]

By (\ref{inverse}), $f_k^N$ is in $C^{4,\alpha}$. By differentiating the above equation twice, we obtain 
\[
\begin{split}
&\nabla^2 f_k^N(\xi,\xi)-\nabla^2\log\rho_{k-1}^N(\xi,\xi)-\frac{K}{N}\tr\left(\left(I-\frac{K}{N}\nabla^2f_k^N\right)^{-1} \nabla^2(\nabla^2f_k^N(\xi,\xi))\right)\\
&-\frac{K^2}{N^2}\tr\left(\left(I-\frac{K}{N}\nabla^2f_k^N\right)^{-1} \nabla^2(\nabla f_k^N(\xi))\left(I-\frac{K}{N}\nabla^2f_k^N\right)^{-1} \nabla^2(\nabla f_k^N(\xi))\right)\\
&=\frac{K}{2N}\nabla^2 |\nabla f_k^N|^2(\xi,\xi) =\frac{K}{N}\left(|\nabla^2f_k^N(\xi)|^2+\left<\nabla(\nabla^2f_k^N(\xi,\xi)), \nabla f_k^N\right>\right). 
\end{split}
\]

By using the same arguments in the proof of (\ref{pos}), we obtain $I-\frac{K}{N}\nabla^2f_k^N\geq 0$. Therefore, if $(x,\xi)$ is a point where $\left<\nabla^2 f_k^N(x)(\xi),\xi\right>$ achieves its minimum $\lambda_k^N$ among all points in $\ntorus\times \mathbb S^n$, then 
\[
\begin{split}
&\frac{K}{N}(\lambda_k^N)^2-\lambda_k^N+a_{k-1}^N\leq 0.  
\end{split}
\]
where $a_k^N=\inf_{(x,\xi)\in\ntorus\times\mathbb S^n}\nabla^2\log\rho_k^N(\xi,\xi)$. 

It follows from this and (\ref{fest}) that 
\[
\begin{split}
&\frac{1-\sqrt{1-\frac{4K}{N}a_{k-1}^N}}{\frac{2K}{N}}\leq \lambda_k^N\\
&\leq \left<\nabla^2 f_k^N(v),v\right>=\left<\nabla^2\log\rho_k^N\left(I+\frac{K}{N}\nabla^2\log\rho_k^N\right)^{-1}v,v\right>
\end{split}
\]
for any unit vector $v$. 

Therefore, by letting $(x,v)$ be a minimizer of $\nabla^2\log\rho_k^N(\xi,\xi)$ in $\ntorus\times\mathbb S^n$, we obtain 
\begin{equation}\label{lambda}
\begin{split}
&\frac{1-\sqrt{1-\frac{4K}{N}a_{k-1}^N}}{2} \leq \frac{\frac{K}{N}a_k^N}{1+\frac{K}{N}a_k^N}. 
\end{split}
\end{equation}

We claim that there is an integer $N_0>0$ depending on $\rho_0$ and a constant $\frac{1}{2}\leq C\leq 1$  such that 
\begin{equation}\label{Har}
a_k^N\geq -\frac{C}{\frac{K}{N}(k+1)}
\end{equation}
for all $N\geq N_0$. 

The constant $C$ will be chosen such that $\frac{(1-C)^2}{2C-1}\leq 1$. The inequality (\ref{Har}) is clearly true for $k=0$ if we choose $N_0$ large enough. Assume that $a_{k-1}^N\geq -\frac{C}{\frac{K}{N}k}$, we show (\ref{Har}) holds. Indeed, we have
\[
\begin{split}
&\frac{\frac{K}{N}a_k^N}{1+\frac{K}{N}a_k^N}\geq \frac{1-\sqrt{1+\frac{4C}{k}}}{2}. 
\end{split}
\]
By (\ref{pos}), $1+\frac{K}{N}a_k^N>0$. Therefore, the above equation is equivalent to
\[
\frac{K}{N}a_k^N\geq \frac{1-\sqrt{1+\frac{4C}{k}}}{1+\sqrt{1+\frac{4C}{k}}} . 
\]
An elementary calculation shows that 
\[
\frac{1-\sqrt{1+\frac{4C}{k}}}{1+\sqrt{1+\frac{4C}{k}}}\geq -\frac{C}{k+1}
\]
is equivalent to 
\[
k\geq\frac{(1-C)^2}{2C-1}
\]
if $1\geq C\geq \frac{1}{2}$. 

Since we chose $C$ close enough to 1, the above inequality holds for all $k\geq 1$. This proves the claim and therefore the theorem. 

\smallskip

\section{Proof of Theorem \ref{Harn}}

It follows from (\ref{fest}) and Theorem \ref{diffHarn} that 
\[
\nabla^2f_k^N\geq -\frac{C}{\frac{K}{N}(k+1-C)}I
\]

It follows from this and (\ref{MA}) that 
\begin{equation}\label{rho1}
\begin{split}
\rho_{k-1}^N&=\rho_{k_1}^N(\phi_{k}^N)\det\left(I-\frac{K}{N}\nabla^2f_{k}^N\right) \\
&\leq \rho_{k}^N(\phi_{k}^N)\left(\frac{k+1}{k+1-C}\right)^n. 
\end{split}
\end{equation}

On the other hand, by (\ref{equal}) and (\ref{ineq}),  
\[
\rho_k^N(\phi_k^N(x))\leq \exp(f_k^N(x))\leq \rho_k^N(y)\exp\left(\frac{N}{2K}d^2(x, y)\right). 
\]
By combining this with (\ref{rho1}), we obtain 
\begin{equation}\label{rho2}
\begin{split}
\rho_{k-1}^N(x)&\leq \left(\frac{k+1}{k+1-C}\right)^n\exp\left(\frac{N}{2K}d^2(x, y)\right)\rho_k^N(y). 
\end{split}
\end{equation}

Let $\gamma:[0,1]\to\ntorus$ be a length minimizing geodesic connecting $x$ and $y$. Let $h=\frac{1}{k_2-k_1+1}$. It follows that 
\[
\begin{split}
\frac{1}{2}d^2(x,y)&=\frac{1}{2}\sum_{m=1}^{k_2-k_1+1}\int_{(m-1)h}^{mh}|\dot\gamma(t)|^2dt\\
&=\frac{1}{2h}\sum_{m=1}^{k_2-k_1+1}d^2(\gamma((m-1)h),\gamma(mh)). 
\end{split}
\]

 It follows from this and (\ref{rho2}) that 
\[
\begin{split}
&\rho_{k_1-1}^N(x)\leq \left(\frac{k_1+1}{k_1+1-C}\right)^n\exp\left(\frac{N}{2K}d^2(x, \gamma\left(h\right))\right)\rho_{k_1}^N(\gamma(h))\\
&\leq \left(\frac{k_1+1}{k_1+1-C}\right)^n...\left(\frac{k_2+1}{k_2+1-C}\right)^n\exp\left(\frac{Nh}{2K}d^2(x, y)\right)\rho_{k_2}^N(y)\\
&\leq \left(\frac{k_2+1}{k_1+1-C}\right)^n\exp\left(\frac{N}{2(k_2-k_1+1)K}d^2(x, y)\right)\rho_{k_2}^N(y). 
\end{split}
\]

Therefore, if $t_1$ and $t_2$ are contained in $\left[\frac{(k_1-1)K}{N}, \frac{k_1K}{N}\right)$ and $\left[\frac{k_2K}{N}, \frac{(k_2+1)K}{N}\right)$, respectively, then 
\[
\begin{split}
&u_{t_1}^N(x)\leq \left(\frac{t_2+\frac{K}{N}}{t_1}\right)^n\exp\left(\frac{d^2(x, y)}{2\left(t_2-t_1-\frac{K}{N}\right)}\right)u_{t_2}^N(y). 
\end{split}
\]

\end{document}